\documentclass{amsart}[12pt]

\usepackage{amssymb,latexsym} 
\usepackage{palatino}
\usepackage{mathrsfs} 

\newcommand{\DynkinF}[5]{
\setlength{\unitlength}{#1 true in}

\begin{picture}(3.6, 1)

\put(.2,.1){\line(1,0){1}}
\multiput(.2,.1)(1,0){4}{\circle*{.2}}
\put(2.2,.1){\line(1,0){1}}
\put(1.2,.06){\line(1,0){1}}
\put(1.2,.14){\line(1,0){1}}
\put(1.43,-.05){$>$}
\put(.08,.27){$#2$}
\put(1.08,.27){$#3$}
\put(2.08,.27){$#4$}
\put(3.08,.27){$#5$}
\end{picture}}

\newcommand{\weightDynkinF}[5]{
\setlength{\unitlength}{#1 true in}

\begin{picture}(3.6, 1)

\put(0.08,0){$#2$}
\put(1.08,0){$#3$}
\put(2.08,0){$#4$}
\put(3.08,0){$#5$}
\end{picture}}

\newcommand{\DynkinEsix}[1]{
\setlength{\unitlength}{#1 true in}

\begin{picture}(4.6, 1)(0,-.55)
\multiput(.2,.1)(1,0){5}{\circle*{.2}}
\put(.2,.1){\line(1,0){4}}
\put(2.2,-.5){\circle*{.2}}
\put(2.2,.1){\line(0,-1){.5}}
\end{picture}}

\newcommand{\weightDynkinEsix}[7]{
\setlength{\unitlength}{#1 true in}

\begin{picture}(4.7, .9)
\put(0.08,0){$#2$}
\put(1.08,0){$#3$}
\put(2.08,0){$#4$}
\put(3.08,0){$#5$}
\put(4.08,0){$#6$}
\put(2.08,-.8){$#7$}

\end{picture}}

\newcommand{\DynkinEseven}[1]{
\setlength{\unitlength}{#1 true in}

\begin{picture}(5.7, 1)(0,-.55)
\multiput(.2,.1)(1,0){6}{\circle*{.2}}
\put(.2,.1){\line(1,0){5}}
\put(2.2,-.5){\circle*{.2}}
\put(2.2,.1){\line(0,-1){.5}}
\end{picture}}

\newcommand{\weightDynkinEseven}[8]{
\setlength{\unitlength}{#1 true in}

\begin{picture}(5.8, .9)
\put(0.08,0){$#2$}
\put(1.08,0){$#3$}
\put(2.08,0){$#4$}
\put(3.08,0){$#5$}
\put(4.08,0){$#6$}
\put(5.08,0){$#7$}
\put(2.08,-.8){$#8$}

\end{picture}}

\newcommand{\DynkinEeight}[1]{
\setlength{\unitlength}{#1 true in}

\begin{picture}(6.7, 1)(0,-.55)
\multiput(.2,.1)(1,0){7}{\circle*{.2}}
\put(.2,.1){\line(1,0){6}}
\put(2.2,-.5){\circle*{.2}}
\put(2.2,.1){\line(0,-1){.5}}
\end{picture}}

\newcommand{\weightDynkinEeight}[9]{
\setlength{\unitlength}{#1 true in}

\begin{picture}(6.8, .9)
\put(0.07,0){$#2$}
\put(1.07,0){$#3$}
\put(2.07,0){$#4$}
\put(3.07,0){$#5$}
\put(4.07,0){$#6$}
\put(5.07,0){$#7$}
\put(6.07,0){$#8$}
\put(2.07,-.8){$#9$}

\end{picture}}

\newcommand{\Z}{{\mathbb Z}}

\newcommand{\orb}{{\mathcal O}}

\DeclareMathOperator{\ad}{ad}

\edef\det{\det\nolimits}

\newcommand{\g}{{\mathfrak g}}
\newcommand{\tor}{{\mathfrak t}}

\newcommand{\p}{{\mathfrak p}}
\newcommand{\levi}{{\mathfrak l}}
\newcommand{\uni}{{\mathfrak u}}

\newcommand{\orbit}{{\mathcal O}}
\newcommand{\orbco}{\tilde{\mathcal O}}

\newcommand{\ph}{\Tilde{\Pi}}
\newcommand{\al}{{\alpha}}

\newcommand{\ro}{{\Phi}}

\newcommand{\lat}{{{\mathcal L}}}

\newcommand{\complex}{{\mathbf C}}

\newcommand{\zz}{\mathbb Z}

\newcommand{\tr}{\mbox{tr}}

\newtheorem{theorem}{Theorem}[section]

\newtheorem{proposition}[theorem]{Proposition}

\newtheorem{conjecture}[theorem]{Conjecture}

\theoremstyle{definition}

\theoremstyle{remark}
\newtheorem{remark}[theorem]{Remark}

\newtheorem{example}[theorem]{Example}

\IfFileExists{srcltx.sty}{\usepackage{srcltx}}

\numberwithin{equation}{section}

\usepackage[latin1]{inputenc}
\usepackage{xspace,amssymb,amsfonts,euscript}
\usepackage{amsthm,amsmath}
\usepackage{euscript}
\input xy \xyoption {all}

\RequirePackage{color}
\definecolor{myred}{rgb}{0.75,0,0}
\definecolor{mygreen}{rgb}{0,0.5,0}
\definecolor{myblue}{rgb}{0,0,0.65}

\usepackage{caption}
\usepackage{subcaption} 

\usepackage{caption}
\usepackage{subcaption} 
\RequirePackage{ifpdf}
\ifpdf
  \IfFileExists{pdfsync.sty}{\RequirePackage{pdfsync}}{}
  \RequirePackage[pdftex,
   colorlinks = true,
   urlcolor = myblue, 
   citecolor = mygreen, 
   linkcolor = myred, 
   pagebackref,
   bookmarksopen=true]{hyperref}
\else
  \RequirePackage[hypertex]{hyperref}
\fi
\RequirePackage{ae, aecompl, aeguill} 

\begin{document}

\title {Irreducible local systems on nilpotent orbits}

\author{Eric N. Sommers} \address{University of
Massachusetts Amherst\\ Amherst, MA 01003}
\email{esommers@math.umass.edu}

\begin{abstract}
Let $G$ be a simple, simply-connected algebraic group over the complex numbers with Lie algebra $\g$.
The main result of this article is a proof that each irreducible representation
of the fundamental group of the orbit $\orb$ through a nilpotent element $e \in \g$ 
lifts to a representation of a Jacobson-Morozov parabolic subgroup of $G$ associated to $e$.
This result was shown in some cases by Barbasch and Vogan in their study of unipotent representations for complex groups and, in general, in an unpublished part of the author's doctoral thesis.

In the last section of the article, we state two applications of this result, whose details will appear elsewhere:
to answering a question of Lusztig regarding special pieces in the exceptional groups (joint work with Fu, Juteau, and Levy); and to computing the $G$-module structure of the sections of an irreducible local system on $\orb$.  A key aspect of the latter application is some new cohomological statements that generalize those in earlier work of the author.
\end{abstract}

\date{December 5, 2016}

\maketitle

\begin{center}
{\it To George Lusztig, on the occasion of his 70th birthday}
\end{center}

\section{Lifting Result} 

Let $G$ be a simple, simply-connected 
algebraic group defined over the complex numbers $\complex$ with
Lie algebra $\g$.
Let $\orbit$ be a nilpotent orbit in $\g$.  Picking $e \in \orbit$,
we can identify $\orbit$ with $G/ G_e$, where $G_e = Z_G(e)$ is
the centralizer in $G$ of $e$ under the adjoint action. 


\subsection{Statement of the result} \label{intro}
This paper is concerned with the 
irreducible representations 
$$\pi:G_e \to GL(V)$$
such that the identity component $G^{\circ}_e$ of $G_e$ is in the kernel of $\pi$.
In other words, $\pi$ descends to an irreducible representation
of the component group $A(e):=G_e / G_e^{\circ}.$
Since this group identifies with the fundamental group of $\orbit$ at
the base point $e$ (in the analytic topology), 
we refer to $\pi$ or its associated bundle 
$G \times^{G_e} V$
over $\orbit$ as an irreducible 
local system on $\orbit$.



Put $e$ in an $\mathfrak{sl}_2$-triple $\{e, h, f \}$ and write
$$\g = \displaystyle\bigoplus_{i \in \Z} \g_i$$
where $\g_i$ is the $i$-eigenspace of $\ad(h)$ on $\g$.  Let 
$\p = \oplus_{i \geq 0} \g_i$.  Let $P$ be the subgroup of $G$ with Lie algebra $\p$ and let $L$ be
the subgroup of $P$ with Lie algebra $\g_0$.  
It is known that $G_e \subset P$ (see \cite{Carter}).
The main result of the paper is 

\begin{theorem} \label{lifting}
Let $(\pi, V)$ be an irreducible representation of $G_e$, trivial on $G_e^{\circ}$.  
Then there exists a representation $(\tilde{\pi}, V)$ of $P$
such that
$\tilde{\pi}|_{G_e} = \pi.$
\end{theorem}

Since $\pi$ is irreducible, any lifting $\tilde{\pi}$ to $P$ 
must also be irreducible.  Therefore the
unipotent radical $U_P$ of $P$ acts trivially on $V$ since the $U_P$-invariants
are stable under the action of $P$ and must be nonzero.
Thus $L$ must act irreducibly on $V$, 
and hence if $\tilde{\pi}$ exists, it can be specified by a highest weight representation $\lambda$ of $L$,
after choosing a maximal torus $T$ in $L$ with $h \in  \tor$
and a Borel subgroup $B$ with $T \subset B \subset P$, where $\tor$ is the Lie algebra of $T$.  

In the classical groups, our proof is by direct construction using partitions;
in the exceptional groups, we use the explicit knowledge of the structure
of $A(e)$ from \cite{sommers:b-c}.  
For applications it is also useful to find $\lambda \in\tor^*$ of minimal length,
subject to a fixed $W$-invariant form on $\tor^*$.

For some orbits, Theorem \ref{lifting}
was proved by Barbasch and Vogan in \cite[\S 9]{bv:unip} as part of their study of unitary representations for complex Lie groups. 
For the cases they consider there is a unified---and somewhat mysterious---explanation for the highest weights that arise. 

Theorem \ref{lifting} in full generality was proved in my PhD thesis \cite{sommers:thesis} under the direction of George Lusztig.   It was not published as it became part of a manuscript that dealt with the graded $G$-module of regular functions on the universal cover of $\orb$.  That manuscript (circa 2005) was distributed on a limited basis. 
In \S \ref{functions} we discuss some of the results from the manuscript.  In \S \ref{special} we give another application of the theorem 
to a question of Lusztig from \cite{lusztig:notes} concerning special pieces in the exceptional groups, which is joint work in progress with Fu, Juteau, and Levy.  

I thank David Vogan for many helpful discussions about Theorem \ref{lifting} and its applications to rings of functions on orbit covers.  I am deeply grateful to George Lusztig, for suggesting that Theorem \ref{lifting} should hold in general, for his supervision of my thesis, and for his continued encouragement and friendship over the years.

\subsection{Algorithm for lifting}

We give a method to establish Theorem \ref{lifting}, which we carry out 
in the exceptional groups in the next section.  Keep the notation from \S \ref{intro}.
Let $\mathfrak m$ denote the span of the $\mathfrak{sl}_2$-triple $\{e, h, f \}$ through $e$.
Then $G_e = Z_G(\mathfrak m) U_e$ is a Levi decomposition,
where $U_e$ is the unipotent radical of $G_e$. 
Now let $V_{\lambda}$ be a highest weight representation of $L$, viewed also as a representation of $P$ where $U$ acts trivially.
Since $U_e \subset U$ and $Z_G(\mathfrak m) \subset L$, we are reduced to studying 
the restriction of $V_{\lambda}$ to the reductive group $Z_G(\mathfrak m)$. 

The first step in showing Theorem \ref{lifting} is to see whether $G^{\circ}_{e}$ acts trivially on $V_{\lambda}$,
which, by the above, is equivalent to $Z^{\circ}_G(\mathfrak m)$ acting trivially.  
This condition is equivalent to a maximal torus $T'$ of $Z_G(\mathfrak m)$, and hence of $G^{\circ}_{e}$, 
acting trivially on $V_{\lambda}$, which is equivalent to 
each weight of $V_{\lambda}$ containing $T'$ in its kernel.


To check this condition, 
we recall some results of Bala and Carter \cite{Carter} and their generalization in \cite{sommers:b-c},
and we refer to those references for proofs.
Let $W = N_G(T)$ be the Weyl group. 
Let $\ro \subset \tor^*$ be the roots determined by $T$ and let
$\Pi$ be the simple roots determined by $B$. 
Let $-\theta$ be the lowest root of $\ro$ given the choice of $\Pi$.
Set $\ph = \Pi \cup \{ -\theta \}$.
Let $\lat$ be the lattice of characters $X^*(T)$ of $T$.
Given $J \subsetneq \ph$, define
$\lat_J$ to be the $\zz$-lattice in $\lat$ spanned by the elements of $J$.
Set $\ro_J = \lat_J \cap \ro$.
Define the subalgebra
$$\g_J = \tor \oplus \bigoplus_{\al \in \ro_J} \g_{\al}$$
where 
$\g_{\al}$ is the root space for $\al$.
Let $G_J$ denote the connected group in $G$ with Lie algebra $\g_J$.

By Bala-Carter \cite{Carter}, there exists a pair $(\g_J, e_1)$ with $J \subset \Pi$ and $e_1 \in \g_J$ distinguished nilpotent.
Then the identity component of $T^{\lat_J}$ is a maximal torus of $G_{e_1}$.
Put $e_1$ in an ${\mathfrak sl}_2$-triple $\{e_1, h_1, f_1 \}$ in $\g_J$ with $h_1 \in \tor$.
Then there exists
$w \in W$ such that $w(h_1) = h$ and 
it follows that $w(e_1) \in \g_2$ (using $w$ also for any lift to $N_G(T)$) and 
so $w(e_1) = x \cdot e$ for some $x \in L$
since $L$ acts transitively on $G \cdot e \cap \g_2$.  
Now $G^{\circ}_e$ acts trivially on $V_{\lambda}$
if and only if 
$G^{\circ}_{x \cdot e} = x (G^{\circ}_e) x^{-1}$ does. 
Hence it is enough to check that a maximal torus
of $G^{\circ}_{w(e_1)}$ acts trivially on $V_{\lambda}$.
Such a maximal torus is given by the identity component of $w(T^{\lat_J})$.
Now a weight $\mu$ is trivial on the identity component of $w(T^{\lat_J})$
if and only if $n \cdot w^{-1}(\mu)$ is trivial on $T^{\lat_J}$ for some positive integer $n$,
which is equivalent to $n \cdot w^{-1}(\mu) \in \lat_J$.  Therefore,

\begin{proposition} \label{trivial_restrict}
Let $\{e, h, f\}$ be an ${\mathfrak sl}_2$-triple used to define $P$ and $L$ 
with $h \in \tor$.  Let $w \in W$ be such that
$w^{-1}(h)$ is the semisimple part of an ${\mathfrak sl}_2$-triple 
$\{e_1, w^{-1}(h), f_1\}$ with $e_1$ distinguished in $\g_J$ for $J \subset \Pi$.
Then $G^{\circ}_e$ acts trivially on $V_{\lambda}$ if and only if
each weight of $V_{\lambda}$ lies in the rational closure of $w {\lat}_J$ in $X^*(T)$.
\end{proposition}

\begin{example}\label{e6}
Consider the nilpotent orbit of type $D_{4}(a_{1})$ in $E_{6}$.
Let 
$J = \{ \al_2, \al_3, \al_4, \al_5 \}$,
using Bourbaki's notation for labeling simple roots,
so that $\g_J$ has semisimple part of type $D_4$.
Consider the semisimple element $h_1 \in \g_J \cap \tor$ given by
$4 \al^{\vee}_{3} + 6 \al^{\vee}_{4} + 4 \al^{\vee}_{5} + 4 \al^{\vee}_{2}$
in the basis of simple coroots of $E_{6}$.
Then $h_1$ completes to an $\mathfrak{sl_{2}}$-triple in $\g_J$
whose nilpotent elements belong to $D_4(a_1)$.
Applying 
$w = s_{4}s_{3}s_{5}s_{2}s_{4}s_{3}s_{5}s_{1}s_{6}$
where $s_{i}:= s_{\alpha_{i}}$ to $h_1$ yields a dominant element $h$
with values $\{ \alpha_i(h) \}$ that give 
the weighted Dynkin diagram $\begin{smallmatrix}
0 & 0 & 2 & 0 & 0 \\
& & 0 & & &
\end{smallmatrix}$ of $D_{4}(a_{1})$.
Hence the Levi subgroup $L$ has semisimple type $A_{2}+A_{2}+A_{1}$
corresponding to all simple roots of $E_6$ except for $\al_4$.

The highest weight representation $V_{\varpi_{2}}$ of $L$ is two-dimensional with weights 
$\mu_{1} = \varpi_{2}, \mu_{2} = \varpi_{2} - \alpha_{2}$.
In a basis of simple roots $\varpi_2$ equals
$\begin{smallmatrix}
1 & 2 & 3 & 2 & 1 \\
& & 2 & & &
\end{smallmatrix}$.
Applying $w^{-1}$ to $\mu_{1}$ and $\mu_{2}$ yields, respectively,
$\begin{smallmatrix}
0 & 1 & 2 & 1 & 0 \\
& & 1 & & &
\end{smallmatrix}$ and
$\begin{smallmatrix}
0 & 1 & 1 & 1 & 0 \\
& & 1& & &
\end{smallmatrix},$
which both lie in $\lat_{J}$, the weights spanned by the roots of $J$.
Hence  $G^{\circ}_e$ acts trivially on $V_{\varpi_{2}}$.
\end{example}

\begin{remark} \label{minuscule}
All the representations in the proof of Theorem \ref{lifting} turn out to be minuscule for $L$ and so the weights of $V_{\lambda}$ are a single orbit under the action of the Weyl group $W_L$ of $L$.  Consequently only the highest weight 
$\lambda$ in Proposition \ref{trivial_restrict} ends up needing to be checked.
\end{remark}

For the rest of the section, assume $G$ is of adjoint type, so that $\lat$ is generated by $\ro$.
We review some results from \cite{sommers:b-c}.
Define $c_{\al}$ by the equation
$\theta = \sum_{\al \in \Pi} c_{\al} \al$
and set $c_{-\theta} = 1$.
Let 
$$d_J = \text{gcd}(c_{\al})_{\al \in \ph - J}.$$
Since $G$ is adjoint, the torsion subgroup of $\lat / \lat_J$ is cyclic 
of order $d_J$ and is generated by
the image of the element 
$$\tau_J = \frac{1}{d_J} (\sum_{\al \in \ph - J}  c_{\al} \al).$$
The center $Z(G_J$) of $G_J$ equals $T^{\lat_J}$
and so the character group of $Z(G_J)$ is isomorphic
to $\lat / \lat_J$.


Let $C$ be a conjugacy class in $A(e)$. 
Then there exists a pair 
$(\g_J, e_1)$
with $J \subsetneq \ph$ and $e_1 \in \g_J$ a nilpotent element
with the following properties:
$e_1$ is distinguished in  $\g_J$;
$e = g \cdot e_1 $ for some $g \in G$;
and any $y \in Z(G_J)$ whose image generates 
$Z(G_J)/Z^{\circ}(G_J)$ has the property that
the image of $gyg^{-1}$ in $A(e)$ lies in $C$.  The trivial conjugacy class $C$
corresponds to the case where $J \subset \Pi$.


Now given that $V_{\lambda}$ is trivial on $G^{\circ}_e$ and thus descends to a representation
of $A(e)$, we can describe how to compute the character of $V_{\lambda}$ as a representation
of $A(e)$ on a conjugacy class $C \subset A(e)$.
First, let $(\g_J, e_1)$ be a pair as above corresponding to $C$.
As before, put $e_1$ in an ${\mathfrak sl}_2$-triple $\{e_1, h_1, f_1 \}$ in $\g_J$ with $h_1 \in \tor$
and let $w \in W$ be such that $w(h_1) = h$.
Then as before $w(e_1) = x \cdot e$ for some $x \in L$. 
Next, choose $s \in T^{\lat_{J}}$ so that $\tau_{J}(s) = \xi$, a primitive $d_J$-th root of unity.
Then $s \in Z(G_J)$ and its image generates the component group of $Z(G_J)$.
Hence, the image of $gsg^{-1}$ in $A(e)$ lies in $C$, where $g:= x^{-1}w$.

Consequently the trace of an element of $C$ on $V_{\lambda}$
can be computed using the element $gsg^{-1}$ and
this trace equals the trace of $wsw^{-1} \in T$ since $x \in L$.
Let $\mu_1, \mu_2, \dots$ be the (not necessarily distinct) weights of $V_{\lambda}$\footnote{That is, 
weight spaces are viewed as one-dimensional, so that weights can be repeated.  
All the weight spaces are, in fact, one-dimensional in the cases considered, see Remark \ref{minuscule}.}.
Then the desired trace 
is given by $\sum \mu_i(wsw^{-1})$, and this can be computed as follows. 
Since $\mu_i$ is trivial on any torus in $G_{w(e_1)}$ and 
since $T^{\lat_J} \subset G_{e_1}$, it follows that 
$w^{-1}(\mu_i)$ is trivial on the identity component of $T^{\lat_J}$.
Hence the image of $w^{-1}(\mu_i)$ in $\lat / \lat_J$ is torsion and so must be 
an integral multiple $a_i$ of $\tau_J$, and so $\mu_i(wsw^{-1}) = (w^{-1} \mu_i)(s)= \xi^{a_i}$.  Therefore,

\begin{proposition}
Let $\{e, h, f\}$ be as before. 
Let $w \in W$ be such that
$w^{-1}(h)$ is the semisimple part of an ${\mathfrak sl}_2$-triple 
$\{e_1, w^{-1}(h), f_1\}$ 
with $e_1$ distinguished in $\g_J$ for $J \subsetneq \ph$.
Let $\{ \mu_i \}$ be the weights of $T$ on a representation $V_{\lambda}$
of $L$ where $G_e^{\circ}$ acts trivially. 
Then for each $i$ there exists $a_i \in \Z$ with $w^{-1}(\mu_i) = a_i \tau_J$ modulo $\lat_J$
and
the trace on $V_{\lambda}$ 
of any element in the conjugacy class of $C$ parametrized by $(\g_J, e_1)$
equals $\sum_{i} \xi^{a_i}$.
\end{proposition}

\begin{example}\label{e6:2}
We compute the character of $V_{\varpi_{2}}$ from the previous example
on the conjugacy class parametrized by $3A_{2}$ 
(the notation refers to the regular nilpotent element in the subalgebra of type $3A_{2}$).  
Here $\tau_{J} = \alpha_{4}$. 
A semisimple element $h_1$ of an $\mathfrak sl_2$-triple for $e_1$ of type $3A_2$ has weighted diagram in $E_{6}$ 
equal to 
$\begin{smallmatrix}
2 & 2 & -6 & 2 & 2 \\
  &   &  2 &  &
\end{smallmatrix}$.
Then 
$w = s_{4} s_{3} s_{5} s_{2} s_{1} s_{4} s_{6} s_{3} s_{5} s_{2} s_{1} s_{4} s_{6} s_{3} s_{5} s_{2} s_{4}$ 
sends
$h_1$ to $h$.
Then $w^{-1} \mu_{1} \equiv - \tau_{J}$ and $w^{-1} \mu_{2} \equiv \tau_{J}$, modulo $\lat_J$,
and therefore the character value is $\xi^{-1} + \xi= -1$,
where in this case $\xi$ is a primitive third root of unity.
In a similar fashion, the character on the conjugacy class parametrized by $A_{3}+2A_{1}$ is computed to be $0$
and so $V_{\varpi_{2}}$ is the irreducible representation of $A(e) \simeq S_{3}$
of dimension two.  In fact, the representations for this orbit fall under the framework in \cite[\S 9]{bv:unip}.
\end{example}




\section{Liftings in the exceptional groups} \label{reps:exceptional}


\subsection{} We carried out the preceding algorithm for the exceptional groups
to prove Theorem \ref{lifting}.  
Let $A^{\text{simp}}(e)$ denote the component group relative to a simply-connected $G$ and 
let $A^{\text{adj}}(e)$ denote the component group relative to the adjoint group $G/Z(G)$. 
The information for the irreducible representations of $A^{\text{adj}}(e)$ are recorded in \S \ref{tables}.

The algorithm begins by first running through the nonzero nodes of the weighted Dynkin diagram of $e$
to see if the corresponding one-dimensional representation $V_{\varpi}$ of $P$
descends to $A^{\text{simp}}(e)$.  
Remarkably this always happens in types $G_{2}, F_{4}, E_{7}, E_{8}$ (and later we will see it does in types $B, C, D$),
even when $A^{\text{simp}}(e)$ is trivial and even when $\varpi$ is not in the root lattice.
Another way to say this is that $Z_{\g}(\mathfrak{m})$, which we know belongs to $\levi$, actually
belongs to $[\levi , \levi ]$.  In many cases, this is true because
$Z_{\g}(\mathfrak{m})$ is semisimple, but in the remaining cases it seems surprising.
This fact is related to the observation by several researchers that in the above cases 
the number of irreducible components of codimension one of the complement of $P \cdot e$ in $\g_2$ 
is the number of nonzero nodes of the weighted Dynkin diagram.
In type $E_{6}$ the same fact holds whenever $A^{\text{simp}}(e)$ is non-trivial
even when $\varpi$ is not in the root lattice.
Consequently, we omit the trivial representation from our tables since
we get a trivial representation of $A^{\text{simp}}(e)$ for each node of the weighted Dynkin diagram
with nonzero value that does not yield a non-trivial representation of $A^{\text{simp}}(e)$.

We also checked in $E_{6}$ (respectively, $E_{7}$) that if $\varpi$ is a fundamental
weight not in the root lattice 
and $\varpi$ corresponds to a node with non-zero value in the weighted Dynkin diagram,
then $V_{3\varpi}$ (respectively, $V_{2\varpi}$) descends to a trivial representation of
$A^{\text{adj}}(e)$.  This implies that the kernel of the representation 
$V_{\varpi}$, viewed as a representation of $A^{\text{simp}}(e)$, 
is a normal subgroup of $A^{\text{simp}}(e)$ isomorphic to $A^{\text{adj}}(e)$.
Hence $A^{\text{simp}}(e)$ must be a split central extension of $A^{\text{adj}}(e)$.
Furthermore, we observe that whenever  $|A^{\text{simp}}(e)| > |A^{\text{adj}}(e)|$, there is always
one such node in the weighted Dynkin diagram.  This implies that 
$A^{\text{simp}}(e)$ is always a split central extension of $A^{\text{adj}}(e)$ in $E_{6}$
and $E_{7}$, which gives a new proof of the splitting in these cases (we note
that there is no splitting in general in types $B$ and $D$, see \S \ref{simply_conn_classical}).
It follows that by tensoring the one-dimensional representations coming from 
these fundamental weights (which are trivial on $A^{\text{adj}}(e)$)
with the representations that we found for $A^{\text{adj}}(e)$ (which are trivial on the image of the
center of $G$),
all the irreducible representations for $A^{\text{simp}}(e)$ in $E_{6}$ and $E_{7}$ are obtained.  
This completes the proof of Theorem \ref{lifting} in the exceptional groups.

\subsection{Tables}  \label{tables}
In the following tables we list weights that yield nontrivial irreducible representations of $A^{\text{adj}}(e)$
in the exceptional groups.  From the discussion in the previous section, 
this allows us to find all weights of $L$ that give rise to irreducible representations of $A^{\text{simp}}(e)$.   
From there, it is easy enough to construct lifts of minimal length.  
The numbering of simple roots is given in types $E$ by
$\begin{smallmatrix}
1 & 3 & 4 & \dots & n \\
  &   &  2 &  &
\end{smallmatrix}$ as in Bourbaki (but different from that in \cite{sommers:thesis}).
We first list the cases where $A^{\text{adj}}(e) \simeq S_2$.

\begin{center}
\begin{tabular}{|c|c|c|} \hline
\multicolumn{3}{|c|}{$F_{4}$} \\ \hline
\multicolumn{1}{|r|}{$\DynkinF{.2}{}{}{}{}$} 
& \multicolumn{1}{|c|}{Bala-Carter}
& \multicolumn{1}{|c|}{Sign rep} 
\\ \hline
\weightDynkinF{.2}{0}{0}{0}{1} & $\Tilde A_{1}$ & $\varpi_{4}$  \\ \hline 
\weightDynkinF{.2}{2}{0}{0}{0} & $A_{2}$ & $\varpi_{1}$  \\ \hline 
\weightDynkinF{.2}{2}{0}{0}{1} & $B_{2}$ & $\varpi_{4}$  \\ \hline 
\weightDynkinF{.2}{1}{0}{1}{0} & $C_{3}(a_{1})$ & $\varpi_{3}$ \\ \hline 
\weightDynkinF{.2}{0}{2}{0}{2} & $F_{4}(a_{2})$ & $\varpi_{2}$ \\ \hline 
\weightDynkinF{.2}{2}{2}{0}{2} & $F_{4}(a_{1})$ & $\varpi_{4}$ \\ \hline 
\end{tabular}
\end{center}

\bigskip

\begin{center}
\begin{tabular}{|c|c|c|} \hline
\multicolumn{3}{|c|}{$E_{6}$} \\ \hline
\multicolumn{1}{|c|}{\DynkinEsix{.18}} 
& \multicolumn{1}{|c|}{Bala-Carter}
& \multicolumn{1}{|c|}{Sign rep}  \\ \hline
\weightDynkinEsix{.18}{0}{0}{0}{0}{0}{2} & $A_{2}$ & $\varpi_{2}$ \\ 
    & &           \\ \hline

\weightDynkinEsix{.18}{2}{0}{2}{0}{2}{0} &  $E_{6}(a_{3})$ & $\varpi_{4}$  \\ 
    & &    \\ \hline

\end{tabular}
\end{center}

\bigskip

\begin{center}
\begin{tabular}{|c|c|c|} \hline
\multicolumn{3}{|c|}{$E_{7}$} \\ \hline
\multicolumn{1}{|c|}{\DynkinEseven{.18}} 
& \multicolumn{1}{|c|}{Bala-Carter}
& \multicolumn{1}{|c|}{Sign rep}
 \\ \hline
\weightDynkinEseven{.18}{2}{0}{0}{0}{0}{0}{0} & $A_{2}$ & $\varpi_{1}$ \\ 
	& & \\ \hline 
\weightDynkinEseven{.18}{1}{0}{0}{0}{1}{0}{0} & $A_{2}+A_{1}$ & $\varpi_{1}$, $\varpi_{6}$   \\ 
	& &  \\ \hline 
\weightDynkinEseven{.18}{0}{1}{0}{0}{0}{1}{1} & $D_{4}(a_{1})+A_{1} $ 
	& $\varpi_{3}$, $\varpi_{2}$-$\varpi_{7}$ \\ 
& &  \\ \hline 
\weightDynkinEseven{.18}{0}{0}{1}{0}{1}{0}{0} & $A_{3}+A_{2} $ & $\varpi_{4}$, $\varpi_{6}$ \\ 
& &  \\ \hline 
\weightDynkinEseven{.18}{2}{0}{0}{0}{2}{0}{0} & $A_{4}$ & $\varpi_{1}$ \\ 
& &  \\ \hline 
\weightDynkinEseven{.18}{1}{0}{1}{0}{1}{0}{0} & $A_{4} + A_{1}$ 
		& $\varpi_{1}$, $\varpi_{4}$ \\ 
& &  \\ \hline 
\weightDynkinEseven{.18}{2}{0}{1}{0}{1}{0}{0} & $D_{5}(a_{1})$ & $\varpi_{4}$, $\varpi_{6}$ \\ 
& &  \\ \hline 
\weightDynkinEseven{.18}{0}{2}{0}{0}{2}{0}{0} & $E_{6}(a_{3})$ & $\varpi_{3}$    \\ 
& &  \\ \hline 
\weightDynkinEseven{.18}{2}{0}{2}{0}{0}{2}{0} & $E_{7}(a_{4})$ & $\varpi_{4}$ \\ 
& &  \\ \hline 
\weightDynkinEseven{.18}{2}{0}{2}{0}{2}{0}{0} & $E_{6}(a_{1})$ & $\varpi_{1}$  \\ 
& &  \\ \hline 
\weightDynkinEseven{.18}{2}{0}{2}{0}{2}{2}{0} & $E_{7}(a_{3})$ & $\varpi_{4}$ \\ 
& & \\ \hline 
\end{tabular}
\end{center}

\bigskip

\begin{center}
\begin{tabular}{|c|c|c|} \hline
\multicolumn{3}{|c|}{$E_{8}$} \\ \hline
\multicolumn{1}{|c|}{\DynkinEeight{.18}} 
& \multicolumn{1}{|c|}{Bala-Carter}
& \multicolumn{1}{|c|}{Sign rep} 
\\ \hline
\weightDynkinEeight{.18}{0}{0}{0}{0}{0}{0}{2}{0} & $A_{2}$ & $\varpi_{8}$\\ & & \\ \hline
                 
\weightDynkinEeight{.18}{1}{0}{0}{0}{0}{0}{1}{0} & $A_{2}+A_{1}$ 
	& $\varpi_{1}$, $\varpi_{8}$\\  & & \\\hline
                 
\weightDynkinEeight{.18}{2}{0}{0}{0}{0}{0}{0}{0} & $2A_{2}$ & $\varpi_{1}$\\  & & \\ \hline
                
\weightDynkinEeight{.18}{1}{0}{0}{0}{1}{0}{0}{0} & $A_{3}+A_{2} $ 
	& $\varpi_{1}$, $\varpi_{6}$\\ & & \\ \hline
                
\weightDynkinEeight{.18}{2}{0}{0}{0}{0}{0}{2}{0} &  $A_{4}$ & $\varpi_{8}$\\ & & \\ \hline
                
\weightDynkinEeight{.18}{0}{0}{0}{0}{0}{0}{0}{2} & $D_{4}(a_{1}) + A_{2}$ 
	& $\varpi_{2}$\\ & & \\ \hline
                
\weightDynkinEeight{.18}{1}{0}{0}{0}{1}{0}{1}{0} & $A_{4} + A_{1}$ 
	& $\varpi_{6}$, $\varpi_{8}$\\ & & \\ \hline
                
\weightDynkinEeight{.18}{1}{0}{0}{0}{1}{0}{2}{0} &  $D_{5}(a_{1})$ & $\varpi_{1}$, $\varpi_{6}$\\  
	& & \\ \hline
                 
\weightDynkinEeight{.18}{0}{0}{1}{0}{0}{0}{1}{0} &  $A_{4} + 2A_{1}$ & $\varpi_{4}$, $\varpi_{8}$\\ 
 & & \\ \hline
                
\weightDynkinEeight{.18}{0}{0}{0}{0}{0}{0}{2}{2} & $D_{4}+A_{2}$ & $\varpi_{2}$\\
 & & \\ \hline
                
\weightDynkinEeight{.18}{2}{0}{0}{0}{0}{2}{0}{0} &  $E_{6}(a_{3})$ & $\varpi_{7}$\\  & & \\ \hline
                
\weightDynkinEeight{.18}{0}{1}{0}{0}{0}{1}{0}{1} &  $D_{6}(a_{2})$ & $\varpi_{3}$, $\varpi_{7}$\\
& & \\  \hline
                
\weightDynkinEeight{.18}{1}{0}{0}{1}{0}{1}{0}{0} & $E_{6}(a_{3})+A_{1}$ & $\varpi_{5}$, $\varpi_{7}$\\ & & \\ \hline
                
\end{tabular}
\end{center}

\newpage
\begin{center}
\begin{tabular}{|c|c|c|} \hline
\multicolumn{3}{|c|}{$E_{8}$} \\ \hline
\multicolumn{1}{|c|}{\DynkinEeight{.18}} 
& \multicolumn{1}{|c|}{Bala-Carter}
& \multicolumn{1}{|c|}{Sign rep} \\ \hline

\weightDynkinEeight{.18}{0}{1}{0}{0}{0}{1}{2}{1}  & $D_{6}(a_{1})$ & $\varpi_{2}$, $\varpi_{3}$\\ 
                 & & \\ \hline
\weightDynkinEeight{.18}{0}{0}{1}{0}{1}{0}{2}{0}  & $E_{7}(a_{4})$ & $\varpi_{4}$, $\varpi_{6}$\\ 
                 & & \\ \hline
\weightDynkinEeight{.18}{2}{0}{0}{0}{2}{0}{2}{0} & $E_{6}(a_{1})$ & $\varpi_{8}$\\ 
                 & & \\ \hline
\weightDynkinEeight{.18}{0}{0}{0}{2}{0}{0}{2}{0} & $D_{5}+A_{2}$ & $\varpi_{5}$\\
                 & & \\ \hline
\weightDynkinEeight{.18}{1}{0}{1}{0}{1}{0}{1}{0}  & $D_{7}(a_{2})$ & $\varpi_{1}$, $\varpi_{4}$\\ 
                 & & \\ \hline
\weightDynkinEeight{.18}{1}{0}{1}{0}{1}{0}{2}{0} & $E_{6}(a_{1})+A_{1}$ & $\varpi_{8}$\\ 
                 & & \\ \hline
\weightDynkinEeight{.18}{2}{0}{1}{0}{1}{0}{2}{0} & $E_{7}(a_{3})$ & $\varpi_{4}$, $\varpi_{6}$\\ 
                   & & \\ \hline              
\weightDynkinEeight{.18}{2}{0}{0}{2}{0}{0}{2}{0} & $D_{7}(a_{1})$ & $\varpi_{5}$\\ 
                 & & \\ \hline
\weightDynkinEeight{.18}{2}{0}{2}{0}{0}{2}{0}{0} & $E_{8}(a_{5})$ & $\varpi_{4}$, $\varpi_{7}$\\ 
 & & \\ \hline
\weightDynkinEeight{.18}{2}{0}{2}{0}{0}{2}{2}{0} & $E_{8}(b_{4})$ & $\varpi_{4}$\\ 
 & & \\ \hline
\weightDynkinEeight{.18}{2}{0}{2}{0}{2}{0}{2}{0} & $E_{8}(a_{4})$ & $\varpi_{4}$, $\varpi_{8}$\\ 
 & & \\ \hline
\weightDynkinEeight{.18}{2}{0}{2}{0}{2}{2}{2}{0} & $E_{8}(a_{3})$ & $\varpi_{4}$\\ 
 & & \\ \hline

\end{tabular}
\end{center}

\newpage
\begin{center}
Sign and standard representations of $A(e)$ when $A(e) \simeq S_{3}$
\end{center}

\begin{center}
\begin{tabular}{|c|c|c|c|} \hline
\multicolumn{4}{|c|}{$A(e) \simeq S_{3}$} \\ \hline
\multicolumn{1}{|c|}{Group}
& \multicolumn{1}{|c|}{Orbit} 
& \multicolumn{1}{|c|}{Sign}
& \multicolumn{1}{|c|}{Standard} \\ \hline
$G_{2}$ & $G_{2}(a_{1})$  & $\varpi_{2}$ & $\varpi_{1}$ \\ 
$E_{6}$ & $D_{4}(a_{1})$ & $\varpi_{4}$ & $\varpi_{2}$ \\ 
$E_{7}$ & $D_{4}(a_{1})$ & $\varpi_{3}$ & $\varpi_{1}$ \\ 
$E_{7}$ & $E_{7}(a_{5})$ & $\varpi_{4}$ & $\varpi_{2} - \varpi_{7}$ \\ 
$E_{8}$ & $D_{4}(a_{1})$ & $\varpi_{7}$ & $\varpi_{8}$ \\ 
$E_{8}$ & $D_{4}(a_{1})+A_{1}$ & $\varpi_{2}$, $\varpi_{7}$ & $\varpi_{8}$ \\ 
$E_{8}$ & $E_{7}(a_{5})$ & $\varpi_{4}$, $\varpi_{6}$ & $\varpi_{2}$ \\ 
$E_{8}$ & $E_{8}(b_{6})$ & $\varpi_{4}$, $\varpi_{8}$ & $\varpi_{2}$ \\ 
$E_{8}$ & $E_{8}(a_{6})$ & $\varpi_{4}$, $\varpi_{7}$ & $\varpi_{2}$, $\varpi_{8}$ \\ 
$E_{8}$ & $E_{8}(b_{5})$ & $\varpi_{4}$ & $\varpi_{2}$ \\ \hline
\end{tabular}
\end{center}

\begin{center}
Character tables when $A(e) \simeq S_{4}, S_{5}$
\medskip

\begin{tabular}{|c|c|c|c|c|} \hline
\multicolumn{5}{|c|}{$F_{4}(a_{3})$} \\ \hline
\multicolumn{1}{|c|}{Conjugacy class} 
& \multicolumn{1}{|c|}{$\varpi_{1}$}
& \multicolumn{1}{|c|}{$\varpi_{2}$}
& \multicolumn{1}{|c|}{$\varpi_{3}$}
& \multicolumn{1}{|c|}{$\varpi_{4}$} \\ \hline

$F_{4}(a_{3})$ & 2 & 1 & 3 & 3\\ 
$A_{3} + \Tilde{A_{1}}$ & 0 & -1 & 1 & -1 \\ 
$A_{2} + \Tilde{A_{2}}$  & -1 & 1 & 0 & 0 \\ 
$B_{4}(a_{1})$ & 2 & 1 & -1 & -1 \\ 
$A_{1} + C_{3}(a_{1})$ & 0 & -1 & -1 & 1 \\  \hline

\end{tabular}

\begin{tabular}{|c|c|c|c|c|c|c|} \hline
\multicolumn{7}{|c|}{$E_{8}(a_{7})$} \\ \hline
\multicolumn{1}{|c|}{Conjugacy class} 
& \multicolumn{1}{|c|}{$\varpi_{8}$}
& \multicolumn{1}{|c|}{$\varpi_{6}$}
& \multicolumn{1}{|c|}{$\varpi_{5}$}
& \multicolumn{1}{|c|}{$\varpi_{7}$}
& \multicolumn{1}{|c|}{$\varpi_{1}$}
& \multicolumn{1}{|c|}{$\varpi_{2}$} \\ \hline

$E_{8}(a_{7})$  & 4 & 4 & 1 & 6 & 5 & 5 \\ 
$A_{5}+A_{2}+A_{1}$ & -1 & 1 & -1 & 0 & 1 & -1 \\ 
$2 A_{4}$ & -1 & -1 & 1 & 1 & 0&  0 \\ 
$D_{5}(a_{1}) + A_{3}$ &0 & 0 & -1 & 0 & -1 & 1 \\
$D_{8}(a_{5})$ & 0 & 0 & 1 & -2 & 1 & 1\\
$E_{7}(a_{5})+A_{1}$ & 2 & -2 & -1 & 0 & 1 & -1 \\
$E_{6}(a_{3})+A_{2}$ & 1 & 1 & 1 & 0 & -1 & -1\\ \hline
\end{tabular}
\end{center}

\section{Lifting in the classical groups} \label{liftings_classical}

In this section we prove Theorem \ref{lifting} in the classical groups.  

\subsection{Types $B, C, D$.}  We first handle the case relative to a group that is symplectic or special orthogonal. 
Fix $\epsilon \in \{ 0 , 1 \}$.   In what follows, all congruences are modulo 2.
Let $V$ be a complex vector space of dimension $N$ with a non-degenerate bilinear
form $\phi: V \times V \to \complex$
satisfying $\phi(v, w) = (-1)^{\epsilon} \phi(w, v)$
for $v, w \in V$.
Let $H$ be the subgroup of $GL(V)$ preserving the form $\phi$,
and let $G$ be its connected component (so $G$ is not assumed to be simply-connected in this section).
Set $n = \lfloor N/2 \rfloor$.
Then $G$ is of type $B_{n}$ when $\epsilon = 0$, $N$ odd;
of type $D_n$ when $\epsilon = 0$, $N$ even;
and of type $C_n$ when $\epsilon = 1$, $N$ even. 

Let $e$ be a nilpotent element in the Lie algebra of one of these groups with corresponding partition 
$\lambda := (\lambda_1 \geq \lambda_2 \geq \dots )$ of $N$.  
We first recall a description of a basis of the component groups 
$A'(e):= Z_H(e)/ Z_H^{\circ}(e)$
and 
$A(e):= Z_G(e)/ Z_G^{\circ}(e).$ 
These are elementary abelian $2$-groups with a natural injective map of $A(e)$ into $A'(e)$.

Let  $$\mathcal B(\lambda):= \{ j \in \mathbf{N} \ | \ \lambda_j > \lambda_{j+1}
\text{ and } \lambda_j \not\equiv \epsilon \}.$$
Fix a normalized basis of $V$ with respect to $e$ as discussed in \cite{hesselink:polarizations}.
The basis consists of $v_{i,j}$ for $1 \leq i \leq \lambda_j$ 
and the action of $e$ is given by 
$e.v_{i,j} = v_{i-1, j}$ for $i >1$ and $e.v_{1,j} = 0$
For $k \in \mathcal B(\lambda)$, 
define an element $b_k \in H$ by 
$b_k.v_{i,j} = v_{i,j} \text{ when } j \neq k$
and
$b_k.v_{i,k} = - v_{i,k}.$
Then $b_k \in Z_H(e)$
and the images of all the $b_k$ in $A'(e)$ give a basis
of $A'(e)$ over $\mathbf{F}_2$.
If $\epsilon = 1$, then $A(e) = A'(e)$.
If $\epsilon = 0$, then $A(e)$ is the subgroup of $A'(e)$ of index two
given by 
$$\{ \sum a_j b_j \ | \ a_j \in \mathbf{F}_2 \text{ and } \sum a_j = 0 \},$$
where we also use $b_j$ for its image in $A(e)$.

Let $p$ (respectively, $q$) be the largest even (respectively, odd) part of $\lambda$.
Let 
\begin{align*} 
E &= \{ p, p-2, \dots ,4, 2 \} \text{ and } \\
O &= \{ q, q-2, \dots ,5, 3 \}.
\end{align*}
For $s \in E \cup O$
define the subspace $F_s$ 
spanned by the vectors $v_{i,j}$ satisfying
$$\lambda_j + 2 - 2i \geq s.$$
The $F_s$ are isotropic and satisfy
$F_s \subset F_{s'}$ whenever $s' \leq s.$
The subgroup of $G$ which fixes this partial isotropic flag is 
a parabolic subgroup $P$ of $G$.
The dimension of $F_s$ equals $\sigma(s)$ where
\[
\sigma(s) := \sum_{\lambda_j \geq s} \left( \lfloor \tfrac{\lambda_j - s}{2} \rfloor + 1 \right).
\]

Let $T$ be a maximal torus of $G$ such that 
the $v_{i,j}$ with $\lambda_j + 1 - 2i > 0$ are weight vectors for $T$.
Then $T \subset P$.  Choose a Borel subgroup $B$ of $G$ with $T \subset B \subset P$.
Next let  $\gamma(\complex^*) \subset T$ be the one-parameter subgroup of $G$ given by
$$\gamma(z).v_{i,j} = z^{\lambda_j+1-2i} v_{i,j}$$
for $z \in \complex^*$.
Then $L := Z_G(\gamma)$ is a Levi subgroup of $P$, which contains $T$.  Now, 
$\gamma$ defines a cocharacter associated to $e$ and therefore 
$Z_G(e, \gamma)$ defines a Levi subgroup of $Z_G(e)$.
Moreover, $P$ is the Jacobson-Morozov parabolic associated to $e$ and $\gamma$.  

For $s \in E \cup O \cup \{ 1 \}$, define
$L_s$ to be the subgroup of $G$ that
preserves the subspace spanned by the
vectors $$v_{i,j} \text{ with } \lambda_j + 2 - 2i = s$$
and is the identity on all other $v_{i,j}$.
Then when $s>1$, $L_s \simeq GL_{d_s}(\complex)$
where 
$$d_s = \sigma(s) - \sigma(s+1) = \# \{ j \ | \ \lambda_j \geq s, \lambda_j \equiv s \}.$$
On the other hand,
$L_1$ is of the same type as $G$, but acting on a vector space of
dimension $\# \{ j \ | \lambda_j \equiv 1 \}.$
Then 
$$L \simeq \prod_{s \in E \cup O \cup \{1\}} L_s.$$

For $s \in E \cup O$,
define a character $\chi_s$ of $T$ such that for  $t \in T$, we have 
$$\chi_s (t) = \prod_{\lambda_j + 2 -2i = s} t_{i,j},$$
where $t.v_{i,j} = t_{i,j}v_{i,j}$.
This character gives rise to the one-dimensional representation of $L$ (and $P$)
which is trivial on $L_k$ for $k \neq s$ and yields the determinant representation of $L_s \simeq GL_{d_s}(\complex)$.
Let $\Xi$ be the elements of $E \cup O$ such that 
$\sigma(s)$ is not equal to $n$ in type $B_n$ and is not equal to $n$ or $n-1$ in type $D_n$.
Then for $s \in \Xi$, 
we have that $\varpi_{\sigma(s)}$ is the weight of a character of $T$, relative to our choice of $B$.
Moreover, $\chi_{s} = \varpi_{\sigma(s)}$ for $s = \text{max} (E \cup O)$
and $\chi_s = \varpi_{\sigma(s)} - \varpi_{\sigma(s+1)}$, otherwise.
Note that $\chi_s$ is conjugate under $W$ to $\varpi_{d_s}$. 


For $s \in \Xi$ let $\pi_s$ be the projection of $L$ onto $L_s$.  Then 
$\pi_s(Z_H(e, \gamma))$ equals 
\begin{align} \label{factors}
& \prod_{m \equiv s, m \geq s } Sp_{r_m}(\complex) 
\text{ when } s \equiv \epsilon \text{ and }\\ 
& \prod_{m \equiv s, m \geq s} O_{r_m}(\complex) \nonumber
\text{ when } s \not\equiv \epsilon,
\end{align}
where $r_m$ is the number of parts of $\lambda$ equal to $m$.
Each factor on the right above 
sits in $L_s$ as a subgroup preserving the 
vectors $v_{i,j}$ with $\lambda_j+ 2 -2i = s$ 
and $m = \lambda_j$.
It follows that $\pi_s(Z_H^{\circ}(e, \gamma))$ consists of matrices of determinant one,
sitting in $L_s \simeq GL_{d_s}(\complex)$. 
Consequently, each one-dimensional representation of $L_s$,
and hence each $\chi_s$,
when viewed as a representation of $P$ trivial on $U_P$,
is trivial on $Z_H^{\circ}(e) = Z_G^{\circ}(e)$ and thus descends to a representation
of $A(e)$.  

In type $C_n$, we add a part equal to zero as the last part of $\lambda$
and then in all types we define $k_{\text{max}}$ 
to be the largest element of $\mathcal B(\lambda)$ (so $\lambda_{k_{\text{max}}} = 0$ in type $C_n$).
Define $\widetilde{\mathcal B}(\lambda) = {\mathcal B}(\lambda) - \{ k_{\text{max}} \}$. 
For $k \in \widetilde{\mathcal B}(\lambda)$, set
\begin{eqnarray*}
\tilde{b}_k = b_k b_{k'} 
\end{eqnarray*}
where $k'$ is minimal for the property that $k' \in \mathcal B(\lambda)$ and $k < k'$ and 
we set $b_{k_{\text{max}}}$ to be the identity in type $C$.
Then the images of the $\tilde{b}_k$, for $k \in \widetilde{\mathcal B}(\lambda)$, 
give a basis of $A(e)$ over $\mathbf F_2$.

Let 
${\mathcal S} \subset \widetilde{\mathcal B}(\lambda)$
and denote by $\bar\chi_{\mathcal S}$ the
one-dimensional representation of $A(e)$ given
by
\begin{eqnarray}
\bar\chi_{\mathcal S}(\tilde{b}_k) = & -1 \text{ \ if \  } k \in \mathcal S \\
\bar\chi_{\mathcal S}(\tilde{b}_k) = & 1 \text{ \ if \  } k \not\in \mathcal S \nonumber,
\end{eqnarray}
allowing $\tilde{b}_k$ to also stand for its image in $A(e)$.
In this way, the characters of $A(e)$ are parametrized by the subsets of $\widetilde{B}(\lambda)$.

\begin{proposition} \label{classical_lift}
Let $e$ be a nilpotent element in types $B_n, C_n, D_n$ with partition $\lambda$.  
Given $\mathcal S \subset \widetilde{B}(\lambda)$, let 
$$\chi_{\mathcal S} := \sum_{j \in {\mathcal S}} \chi_{\lambda_j}.$$
Then $\chi_{\mathcal S}$ is a lifting to $P$ of 
$\bar\chi_{\mathcal S}$.  This lift is of minimal length.
\end{proposition}

\begin{proof}
By the above discussion, each $\chi_s$ for $s \in \Xi$ descends to 
a character of $A(e)$.  
Next, it is clear from \eqref{factors} that for
$s \equiv \epsilon$ that $\chi_s$ is trivial on $A(e)$.  Also this holds
for $s \leq \lambda_{k_{\text{max}}}$.  For other $s$, we have 
$\chi_s = \chi_{\lambda_k}$ where 
$ \lambda_{k'} < s \leq \lambda_k$ for $k', k \in \mathcal B(\lambda)$.
Moreover for $j \in \widetilde{B}(\lambda)$ it is clear that $\chi_{\lambda_j}$ restricts to $\bar{\chi}_{\{ j \}}$.
This completes the proof except for the claim of minimal length.

In the usual inner produce on weights of $G$, the $\chi_s$ are mutually orthogonal.  They also generate the lattice of one-dimensional characters of $L$, together with possibly one additional weight (in type $B$, when there is exactly one odd $\lambda_j$) and up to two additional weights (in type $D$, when there are at most two odd  $\lambda_j$'s).  These additional weights descend to $A(e)$ and are trivial (this will also follow from results in the next section).  The additional weight(s) can be chosen to be orthogonal to all $\chi_s$ (and to each other).  The orthogonality of these basis elements of the lattice then implies that the 
$\chi_{\mathcal S}$ are minimal lifts.

\end{proof}

\subsection{Liftings for simply-connected type in the classical groups}{ } \label{simply_conn_classical}

It remains to treat the case of the spin groups, which are 2-fold covers of the groups $G$ 
from the previous section in types $B$ and $D$.   Fix $G$ of type $B_n$ ($n \geq 2$) 
or $D_n$ ($n \geq 3$) from the previous section.
Let $\tilde{G}$ denote the spin group covering of $G$ with isogeny $f: \tilde{G} \to G$.  
Let $\tilde{L} = f^{-1}(L)$ and $\tilde{P} = f^{-1}(P)$.
The kernel of $f$ is a group of order $2$ generated by an element $c$ in the center of $\tilde{G}$.
Clearly, $f({\tilde G}_{e}) = G_e$
and $f({\tilde G}^{\circ}_{e}) = G^{\circ}_e$
and so $f$ induces a surjection of 
$$A_{\tilde G}(e):= {\tilde G}_{e}/ {\tilde G}^{\circ}_{e}$$
onto $A_G(e)$ with a kernel
that is either trivial, or non-trivial and generated by the image of $c$.

From \cite{Lusztig:intersect} (see also \cite{sommers:thesis}), 
the kernel is non-trivial if the odd parts of the partition $\lambda$ of $e$ occur with multiplicity one.
To that end, let $e$ be a nilpotent element with partition $\lambda$ such that $r_i \in \{ 0, 1 \}$ when $i$ is odd.
Note that in type $D$ such an element $e$ could be very even. 
Retain the notation of the previous section.
Let $m = |\widetilde{\mathcal B}(\lambda) | + 1$.
Then $|A_{\tilde G}(e)| = 2^m$. 



\begin{proposition}
In type $B_n$ the representation $V_{\varpi_n}$ of $\tilde L$ of highest weight ${\varpi_n}$ has dimension $2^{\frac{m-1}{2}}$.  
When lifted to $\tilde{P}$, it descends to an irreducible representation of $A_{\tilde G}(e)$.
Together with the $2^{m-1}$ characters of $P$ in Proposition \ref{classical_lift}
that descend to $A_G(e)$,
these form a complete set of irreducible representations of $A_{\tilde G}(e)$.
\end{proposition}

\begin{proof}
Let $\ro_J$ be of type 
$$A_{a_1} \times \dots \times A_{a_k} \times B_{q}$$
where $2q+1$ is the sum of the odd parts of $\lambda$
and the even parts of $\lambda$ are
$$[a_1+1, a_1+1, \dots, a_k+1, a_k+1].$$
Assume that the factors of $\ro_J$ sit in $G$ in the order written above.  

Take $e_1 \in \g_J$
so that $e_1$ is regular in the type $A$ factors and has 
partition in the $B_l$ factor consisting of the odd parts of $\lambda$.
Then $e_1$ is distinguished in the Levi subalgebra $\g_J$.
The Dynkin element $h_1$ in $\g_J$ is given by (in the standard basis for coweights)
\begin{eqnarray*}
(a_1, a_1\!-2, \dots, -a_1+2, -a_1, \dots
a_k, a_k\!-2, \dots, -a_k+2, -a_k, 
y_1, \dots, y_{l}, \overbrace{0, \dots , 0}^{\frac{m-1}{2}})
\end{eqnarray*}
where $y_i$ are positive even integers listed in nonincreasing order.
We can conjugate $h_1$ to be dominant using an element $w \in W$ which flips all signs
of the negative elements in $h_1$ and 
then permutes the nonzero elements of $h_1$ into nonincreasing order.

We can now apply Proposition \ref{trivial_restrict} to show that $V_{\varpi_n}$
is trivial on ${\tilde G}^{\circ}_e$.
First observe that $V_{\varpi_n}$ restricts to give the spin representation of $f^{-1}(L_1)$, a spin group of type $B_{\frac{m-1}{2}}$.  Hence the representation is $2^{\frac{m-1}{2}}$-dimensional and all weights are $W_L$-conjugate.  Therefore, it is enough to check 
that the weight $\varpi_n$ is trivial on a maximal torus of ${\tilde G}^{\circ}_e$.
Next, write $\varpi_n$ in the standard basis for weights in $B_n$ as 
$\varpi_n = (\frac{1}{2}, \dots, \frac{1}{2}).$
Then $w^{-1}(\varpi_n)$ equals
\begin{eqnarray*}
(\overbrace{\frac{1}{2}, \dots, \frac{1}{2}}^{\frac{a_1+1}{2}}, 
\overbrace{-\frac{1}{2}, \dots, -\frac{1}{2}}^{\frac{a_1+1}{2}},
\dots, 
\overbrace{\frac{1}{2}, \dots, \frac{1}{2}}^{\frac{a_k+1}{2}}, 
\overbrace{-\frac{1}{2}, \dots, -\frac{1}{2}}^{\frac{a_k+1}{2}},
\overbrace{\frac{1}{2}, \dots, \frac{1}{2}}^{l},
\overbrace{\frac{1}{2}, \dots,\frac{1}{2}}^{\frac{m-1}{2}}).
\end{eqnarray*}
Twice this weight lies ${\mathcal L}_J$.  Hence 
$V_{\varpi_n}$ descends to a representation of $A_{\tilde G}(e)$.

Now the weights of $V_{\varpi_n}$ are not weights of $G$; hence they 
take the value $-1$ on the central element $c$ of $\tilde G$.
So $c$ acts by $-1$ on $V_{\varpi_n}$.  
Now for $x \in A_{\tilde G}(e)$ with $x \neq 1 ,c $, we have that $x$ and $xc$ are conjugate in 
$A_{\tilde G}(e)$ 
by \cite{Lusztig:intersect} (see also \cite[Proposition 26]{sommers:b-c}). 
It follows that 
$$\tr(x, V_{\varpi_n}) = \tr(xc, V_{\varpi_n}) = -\tr(x, V_{\varpi_n})$$
and so 
$\tr(x, V_{\varpi_n}) = 0$
for $x \neq 1 ,c $.
Therefore the inner product of the character of $V_{\varpi_n}$ with itself
is computed just on $1, c \in A_{\tilde G}(e)$,
yielding 
$$\frac{1}{|A_{\tilde G}(e)|}((2^{\frac{m-1}{2}})^2 + (-2^{\frac{m-1}{2}})^2) = \frac{1}{2^m} 2^m = 1,$$
which shows that $V_{\varpi_n}$ is an irreducible representation of $A_{\tilde G}(e)$.
Finally, the sum of the squares of the dimensions of the
distinct irreducible representations of $A_{\tilde G}(e)$ that we have constructed
equals
$$1 \cdot (2^{\frac{m-1}{2}})^2 + 2^{m-1} \cdot (1^2) = 2^{m-1} + 2^{m-1} = 2^{m},$$
which shows that we have found them all.
\end{proof}

To obtain a weight of minimal length for the extra representation, we can subtract off the fundamental weight for the the largest nonzero node $i$ of the weighted Dynkin diagram with $i \leq n-1$,
if such a node exists.

\begin{proposition}
In type $D_n$ when $m \geq 2$, the representations
$V_{\varpi_{n}}$, $V_{\varpi_{n-1}}$ of $\tilde P$ are both of dimension $2^{\frac{m}{2}-1}$
and descend to give irreducible representations of $A_{\tilde G}(e)$.
Together with the $2^{m-1}$ irreducible representations of $P$ in Theorem \ref{classical_lift}
which descend to $A_G(e)$,
these form a complete set of irreducible representations of $A_{\tilde G}(e)$.
\end{proposition}

\begin{proof}
In this case $L_1$ is of type $D_{\frac{m}{2}}$
and $V_{\varpi_{n}}$, $V_{\varpi_{n-1}}$
restrict to give the two half-spin representations of $f^{-1}(L_1)$,
which are of dimension $2^{\frac{m}{2}-1}$.
In case $m=2$, we think of $D_1$ as the center of $\tilde G$
and the half-spin representations are one-dimensional representations which
are non-trivial on $c$.

As in type $B$, the weights of these representations are $W_L$-conjugate and so a similar argument yields
that ${\tilde G}^{\circ}_e$ acts trivially and that $c$ acts by $-1$.
Next, the center of $A_{\tilde G}(e)$ coincides with the image of the center of $\tilde G$ and is of order $4$.
Moreover for $x \in A_{\tilde G}(e)$, with $x$ not in 
the center of $A_{\tilde G}(e)$, the elements $x$ and $xc$ are conjugate in 
$A_{\tilde G}(e)$ \cite{Lusztig:intersect} (see also \cite[\S 4]{sommers:b-c}).  Thus
the character of these representations are zero away from the four
central elements of $A_{\tilde G}(e)$ and the images of central elements of $\tilde G$ act by scalars
since $V_{\varpi_{n-1}}$ and $V_{\varpi_{n}}$ are irreducible representations of $f^{-1}(L)$.
The scalars must be roots of unity and so the inner product of the character of either representation 
with itself is equal to $1$
For the last statement we have the sum of the squares of the known irreducibles is 
$$2 (2^{\frac{m}{2}-1})^2 + 1^2 (2^{m-1}) = 2^{m-1} + 2^{m-1} = 2^{m}$$
as desired.
\end{proof}

For the case where $m=1$, that is, when $e$ is very even, either $V_{\varpi_{n-1}}$ 
or $V_{\varpi_{n}}$ will give the desired one-dimensional representation, depending on which node of the 
weighted Dynkin diagram is nonzero.

To obtain a weight of minimal length for the two extra representations, we can subtract off the fundamental weight for the the largest nonzero node $i$ of the weighted Dynkin diagram with $i \leq n-2$,
if such a node exists.

\subsection{Type $A$}  Let $G$ be the special linear group $SL_{l}(\complex)$.  Let  
$\lambda := (\lambda_1 \geq \lambda_2 \geq \dots \geq \lambda_k)$ be a partition of $l$ 
and let $d$ be the greatest common divisor of $\lambda_1, \lambda_2, \dots, \lambda_k$.  
Let $e$ be a nilpotent element corresponding to $\lambda$. Then $A(e) \simeq \zz / d \zz$.
Let $q = \frac{l}{d}$.  It is easy to check that for $\varpi_{jq}$ fulfill the conditions of Theorem \ref{lifting}
for $0 \leq j \leq d-1$ and that these weights are of minimal length.

%

\section{Applications}

\subsection{Special pieces}\label{special}

We refer to \cite{lusztig:notes} for definitions.  Assume here that $G$ is of adjoint type.
Let $\orb_{\text{sp}} = \orb \cup \orb'$ be a special piece in an exceptional group such that $A(e) \simeq S_2$ for $e$ in the special orbit $\orb$. 
Let $\lambda$ be a minimal lift to $P$ of the sign representation of $A(e)$ as in Theorem \ref{lifting}.
Let $\mu$ be the dominant weight that is $W$-conjugate to $\lambda$ and let $V$ be a representation of $G$ of highest weight $\mu$.
Consider the $G$-orbit $Z$ of $(e, v) \in \g \oplus V$, where $v \in V$ is a weight vector of weight $\lambda$.
Then the stabilizer of $(e,v)$ in $G$ equals $G^{\circ}_e$ and thus $Z \simeq G/ G^{\circ}_e$.   Let $\overline Z$ be the closure
of $Z$ in $\g \oplus V$.

From the tables in \cite{fjls}, one observes that there is a transverse slice $\mathcal S$ in $\overline \orb$ to a point $e' \in \orb'$ that is isomorphic to the closure of the minimal orbit $X_{\text{min}}$ in a symplectic group of type $C_n$.  
Let $p: \overline{Z} \to \overline{\orb}$ be the restriction to $\overline{Z}$ of the projection of $\g \oplus V$ onto the first factor.

In work in progress with Fu, Juteau, and Levy, we checked that 

\begin{proposition} \label{S2_special}
The pre-image $p^{-1}(\mathcal S)$ is isomorphic to $\complex^{2n}$
and $p^{-1}(\mathcal S)$ identifies with the affinization of the universal cover of $X_{\text{min}}$.
\end{proposition}

The proof amounts to checking that when $V$ is restricted to a certain reductive centralizer of $e'$ (which happens to contains a simple factor of type $C_n$), the vector $v$ lies in a subrepresentation isomorphic to the defining representation of $Sp_{2n}(\complex)$.
 A consequence of Proposition \ref{S2_special} is that $p^{-1}(\orb_{\text{sp}} )$ is smooth and satisfies the desired properties of $C^{\dagger}$ from \cite[\S 0.6]{lusztig:notes}.   A statement similar to Proposition \ref{S2_special} is expected to hold for the other special pieces in the exceptional groups.

\subsection{Module structure of sections of local systems} \label{functions}

In \cite{jantzen} Jantzen gives an exposition of results in  McGovern \cite{mcgovern:regfn} and Graham \cite{graham}, results we propose to generalize, and we refer there for definitions.
Let $\orbco:= G / G^{\circ}_e$ denote the universal cover of $\orb$, which carries a left action of $G$ and a right action of $A(e)$.   The $\complex^*$-action on $\orb$ lifts compatibly to 
$\orbco$ so that the regular functions $\complex[\orbco]_{(n)}$ of degree $n$ on $\orbco$ are a $G \times A(e)$-module.  
Let $\hat{A}(e)$ refer to the isomorphism classes of irreducible representations of $A(e)$.
We can write
$$\complex[\orbco]_{(n)} \simeq \bigoplus_{(\pi,V) \in \hat{A}(e)} V \otimes \Gamma^n_V,$$
where  $\Gamma^n_V$ are the degree $n$ sections of the vector bundle 
$G \times^{G_e} V \to \orbit$ corresponding to a lift of $V$ to $G_e$, which is trivial on $G^{\circ}_e$.
We note that the grading is such that $\Gamma^{1}_{\text{trivial}} = 0$ and $\Gamma^{2}_{\text{trivial}}$  is isomorphic to the adjoint representation.

It is expected that Theorem \ref{lifting} can be used to determine  $\Gamma^n_V$ as follows.
Let $(\pi, V) \in \hat{A}(e)$.  Let $(\widetilde{\pi},V)$ be a lifting to $P$ as in 
Theorem \ref{lifting}.   Let $\lambda$ be a lowest weight of $\widetilde{\pi}$ so that $|\lambda|$ is minimal among lowest weights of 
all possible lifts of $V$.  In fact $\lambda$ can always be taken to be a sum of positive roots.   
Let $\complex_{\lambda}$ be the one-dimensional representation of $B$ defined by $\lambda$.  

\begin{conjecture} \label{main_conj}
For $\lambda$ as above, we have
\begin{equation} \label{monty}
\Gamma^n_V \simeq H^0(G/B, S^{(n + \langle \lambda, \chi \rangle)/2 } \g_{\geq 2}^* \otimes \complex_{\lambda})
\end{equation}
for all $n \in \Z$ and 
\begin{equation}  \label{hinich-panyu}
H^i(G /B, S^{n} \g_{\geq 2}^* \otimes \complex_{\lambda}) = 0
\end{equation}
for $i>0$ and for all $n$.
\end{conjecture}
\noindent Here, $\langle \ , \ \rangle$ is the pairing of characters and cocharacters of $T$, and $\chi \in T$ is the cocharacter associated to $e$ and $P$.
It is well-known that in cases where the conjecture holds that there are explicit formulas for the multiplicity of a highest weight representation of $G$ in $\Gamma^n_V$ using 
Lusztig's $q$-analogue of Kostant's weight multiplicity.

We are able to prove a variation of the result, and hence still compute $\Gamma^n_V$, in the following two situations:  (1) when $G$ is of type $A$ and $e$ is any nilpotent element, or 
(2) when $e$ is even and $V$ is one-dimensional.   
In the variation we may have to replace $\g_{\geq 2}^*$ by $\uni_{P'}$, the Lie algebra of the unipotent radical of a parabolic subgroup $P'$ such that $e$ belongs to the Richardson orbit in $\uni_{P'}$.  Moreover in case (2), we may have to replace $\lambda$ by a weight $\mu \in X^*(P')$ which is $W$-conjugate to $\lambda$ and dominant.
The proof uses cohomological statements of the kind proved in \cite{sommers:a_vanishing} and \cite{sommers:normalityD}.
We finish with an example that can be treated using results already in the literature.

\begin{example}
Let $G$ of type $A_5$ and $e$ have partition $(4,2)$.   Consider the weight $\lambda = \varpi_3$ in Theorem \ref{lifting} which corresponds to the non-trivial representation $V$ of $A(e) \simeq \Z/2\Z$.
Now  $2 \varpi_3$ equals $\begin{smallmatrix} 1 & 2 & 3 & 2 & 1 \end{smallmatrix}$ in the basis of simple roots. 
For $S \subset \{1, \dots, |\Pi| \}$, write $\uni_S$ for the Lie algebra of the unipotent radical of the parabolic subgroup containing $B$ 
whose maximal roots are the $-\al_i$ for $i \in S$.
Now 
\begin{align}
& H^i(S^{n} \uni_{1,3,5}^* \otimes -2 \varpi_3) =  H^i(S^{n-2} \uni_{1,3,5}^* \otimes  \begin{smallmatrix} \!-1 & \!-1  & \!-1 & \!-1 & \!-1 \end{smallmatrix}) =  \\
& H^i(S^{n-3} \uni_{2,3,5}^* \otimes  \begin{smallmatrix} 0 & 0  & \!-1 &\! -1 &\! -1 \end{smallmatrix}) =  \nonumber
 H^i(S^{n-4} \uni_{2,3,4}^* \otimes  \begin{smallmatrix} 0 & 0  & -1 & 0 & 0 \end{smallmatrix}) = \\
& H^i(S^{n-5} \uni_{2,3,4}^*) = H^i(S^{n-5} \uni_{1,3,5}^*) \nonumber
\end{align}
for all $i \geq 0$ and all $n$, using the main theorem in \cite{sommers:a_vanishing} five times (we only need this result for $i=0$).
By the proof of Proposition 1.4 in \cite{graham} it follows that $\Gamma^n_V = H^0(S^{(n+5)/2}(\uni_{1,3,5}^*) \otimes - \varpi_3)$.

Similarly, $H^i (S^{(n+5)/2}(\uni_{1,3,5}^*) \otimes - \varpi_3) = H^i (S^{(n-5)/2}(\uni_{1,3,5}^*) \otimes \varpi_3)$ for all $i\geq 0$ and all $n$.  
The latter modules are known to vanish for $i >0$ and all $n$ since $\varpi_3$ is dominant.
Hence in this case the conjecture itself holds.
\end{example}

%
%
%
%
%
%
%
%
%

\bibliography{sommers_vanishing}

\ifx\undefined\bysame
\newcommand{\bysame}{\leavevmode\hbox to3em{\hrulefill}\,}
\fi
\begin{thebibliography}{10}

\bibitem{bv:unip}
D.~Barbasch and D.~Vogan, {\em Unipotent representations of complex semisimple
  Lie groups}, Ann. of Math. {\bf 121} (1985), 41--110.

\bibitem{Carter}
R.~W. Carter, {\em Finite groups of {L}ie type}, Wiley Classics Library, John
  Wiley \& Sons Ltd., Chichester, 1993, Conjugacy classes and complex
  characters, Reprint of the 1985 original, A Wiley-Interscience Publication.

\bibitem{fjls}
B.~Fu, D.~Juteau, P.~Levy, and E.~Sommers, {\em Generic singularities of
  nilpotent orbit closures}, Advances in Mathematics {\bf 305} (2017), 1--77.

\bibitem{graham}
W.~A. Graham, {\em Functions on the universal cover of the principal nilpotent
  orbit}, Invent. Math. {\bf 108} (1992), no.~1, 15--27.

\bibitem{hesselink:polarizations}
W.~H. Hesselink, {\em Polarizations in the classical groups}, Math. Z. {\bf
  160} (1978), no.~3, 217--234.

\bibitem{jantzen}
J.~C. Jantzen, {\em Nilpotent orbits in representation theory}, Lie theory,
  Progr. Math., vol. 228, Birkh\"auser Boston, Boston, MA, 2004, pp.~1--211.

\bibitem{Lusztig:intersect}
G.~Lusztig, {\em Intersection cohomology complexes on a reductive group},
  Invent. Math. {\bf 75} (1984), no.~2, 205--272.

\bibitem{lusztig:notes}
\bysame, {\em Notes on unipotent classes}, Asian J. Math. {\bf 1} (1997),
  194--207.

\bibitem{mcgovern:regfn}
W.~M. McGovern, {\em Rings of regular functions on nilpotent orbits and their
  covers}, Invent. Math. {\bf 97} (1989), no.~1, 209--217.

\bibitem{sommers:thesis}
E.~Sommers, {\em Nilpotent orbits and the affine flag manifold}, Ph.D. thesis,
  Massachusetts Institute of Technology, 1997.

\bibitem{sommers:b-c}
\bysame, {\em A generalization of the Bala-Carter theorem for nilpotent
  orbits}, Internat. Math. Res. Notices (1998), no.~11, 539--562.

\bibitem{sommers:normalityD}
\bysame, {\em Normality of very even nilpotent varieties in {$D_{2l}$}}, Bull.
  London Math. Soc. {\bf 37} (2005), no.~3, 351--360.

\bibitem{sommers:a_vanishing}
\bysame, {\em Cohomology of line bundles on the cotangent bundle of a
  {G}rassmannian}, Proc. Amer. Math. Soc. {\bf 137} (2009), no.~10, 3291--3296.

\end{thebibliography}
\bibliographystyle{pnaplain}

\end{document}